\documentclass [12pt]{article}      
\usepackage{amsmath} 
\usepackage{amsthm}
\usepackage{latexsym}
\usepackage{pb-diagram}
\usepackage{amssymb}
\usepackage{bbm}
\usepackage{color}
\usepackage{graphicx}
\usepackage{hyperref}
\usepackage{cancel}
\hypersetup{
    colorlinks,
    citecolor=blue,
    filecolor=blue,
    linkcolor=blue,
    urlcolor=blue
}
\hypersetup{linktocpage}
\usepackage{enumerate}
\usepackage{todonotes}
\makeatletter
\let\r@eqnnum\@eqnnum
\input{leqno.clo}
\let\l@eqnnum\@eqnnum
\newcommand{\leqnos}{\let\@eqnnum\l@eqnnum}
\newcommand{\reqnos}{\let\@eqnnum\r@eqnnum}
\reqnos
\makeatother

\newtheorem*{theorem*}{Theorem}

\newtheorem*{corollary*}{Corollary}

\newcommand{\la}{\langle}
\newcommand{\ra}{\rangle}

\newcommand{\poZ}{\mathbb Z}
\newcommand{\nn}{{\mathbb N}}

\newcommand{\mcb}{\mathcal B}

\newcommand{\red}[1]{{\textcolor{red}{#1}}}

 \title{Where do you put a telescope? How do you understand Covid concentrations?}
 \author{Matthew Foreman\thanks{Foreman is partially supported by US National Science Foundation grant  DMS-2100367.}}

\begin{document}
 \maketitle
 \noindent{\red{This paper is being edited and shorted considerably for possible publication.}}
 \tableofcontents
 \begin{abstract}
This note discusses dynamical systems--systems that evolve through time.  We start with two contemporary examples illustrating the \emph{qualitative} and the \emph{quantitative} behavior of dynamical systems. These are two broad categories, usually called the study of the \emph{smooth behavior} and \emph{ergodic theory}.  We then introduce the technical framework necessary to state the problems mathematically.  Finally we show that the problems are unsolvable in a rigorous  sense.
\end{abstract}

 \section{Two contemporary examples}
 
\noindent\textbf{Positioning Satellites} Where should you position a satellite so that it stays in a stable orbit? How do you knock an asteroid out of its current trajectory?

 Questions about orbiting bodies arose naturally since Isaac Newton.  Perhaps the most famous is the \emph{n-body problem}. In the late $19^{\mbox{th}}$-century the King of Sweden, inspired by work of Poincar\'{e},  announced a prize for its solution.  It was phrased in the following way (\cite{wikipediapoinc}):

	\begin{quotation}
	\noindent Given a system of arbitrarily many mass points attract each according to Newton's 
	law, 
	under the 
	assumption that no two points ever collide, try to find a representation of the coordinates of 
	each 
	point as a series in a variable that is some known function of time and for all of whose values 
	the
	 series converges uniformly.
	\end{quotation}
A general solution to this problem (even with $n=3$) is impossible (in the sense of having a closed form solution). (See \cite{Saari}.)
So \emph{how do} you put a satellite telescope in orbit so that it stays in place for significantly long time?

In special  cases the three body problem is well approximated by the two body problem which is effectively solvable (\cite{Cornish}). It is 
relevant that when launching a satellite into orbit around the earth and the sun, the mass of the satellite is 
inconsequential relative to the other two masses.

For the earth-sun two body problem there are 5 \emph{Lagrangians}, locations where the relative gravity of the earth and the sun balance out exactly.  So a tiny particle placed exactly on one of those 5 points would stay exactly there. They are \emph{fixed points}.  But  \dots you can't put a satellite in a completely exact location and satellites are not point-masses.

	\begin{figure}[!h]
	\centering
	\includegraphics[height=.25\textheight]{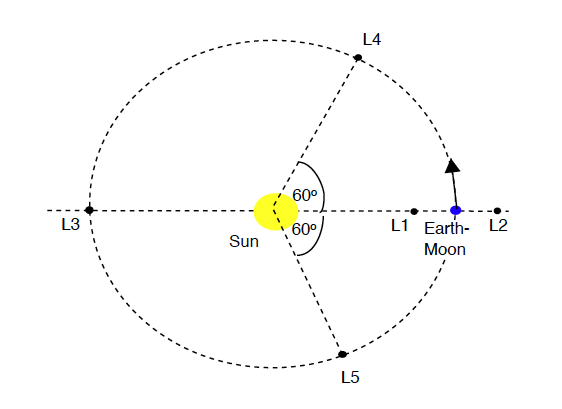}
	\caption{Earth-Sun Lagrangians (\cite{BRw})}
	\label{lagrangians}
	\end{figure}
 Instead you might hope for an \emph{stable point}--a point where small errors get corrected by 
gravity. 
In figure \ref{lagrangians} 
we see that there are 5 earth-sun Lagrangians.  L4 and L5 are stable, but L1, L2 and L3 are not. 

So why not put the satellite L4 or L5?  What if the satellite is carrying the James Webb telescope? For the telescope to function it has to be in the earth's shadow. The only Lagrangian in the earth's shadow is L2 which is unstable.  The strategy is to position the satellite in L2, and use onboard fuel to return the satellite to as close to L2 as possible, thereby approximating stability.

Launched on December 25, 2021, the James Webb telescope was launched into orbit with the target of being as close to L2 as possible. The initial hope was that it would be close enough that the on-board fuel needed to maintain a position near the stable point would be sufficient to keep the telescope operational for 5 years. In fact the launch was so successful at positioning that the current estimate was increased to 20 years. (\cite{SC}). 

	\begin{figure}[!h]
	\centering
	\includegraphics[height=.22\textheight]{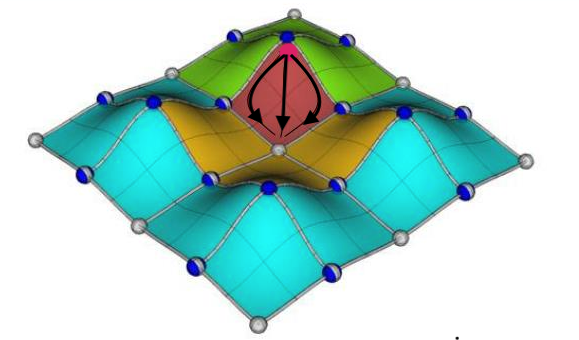}
	\caption{Qualitative behavior of a flow.  \cite{MorseSmale}}
	\label{MorseSmale}
	\end{figure}

\noindent\textbf{Qualitative behavior} Qualitative behavior is exemplified by stable points, unstable points, 
attracting points, points that attract from one direction and repel in another, asymptotic behavior and so 
forth.  Figure \ref{MorseSmale}, showing a Morse-Smale flow,
 illustrates this.  Points travel downward along the surface following the minimal path. (Think of water flowing down a hill.) The local 
minimums are stable attracting points, the local maximums are stable repelling points.  The saddle points 
have some points flowing towards them and an orthogonal collection flowing away.

\noindent\textbf{Quantitative Behavior} Starting in late 2019 and early 2020 the Covid pandemic started. It is spread from individual to individual, but the public health interest is heavily involved with the collective statistical behavior.  Quite striking is the similarity of curves of incidence rates of Covid in many localities.  Is there a mathematical model that predicts the statistics?

The standard model is called the SIRS model, where SIRS stands for Susceptible, Infectious and Recovered.  It has been elaborated on to allow for variable 
time after exposure for an individual to become contagious, leading to the SEIRS model 
(\cite{BSK}, \cite{DP}).
 But \emph{do the two models make the same predictions?}
	\begin{figure}[h]
	\centering
	\includegraphics[height=.25\textheight]{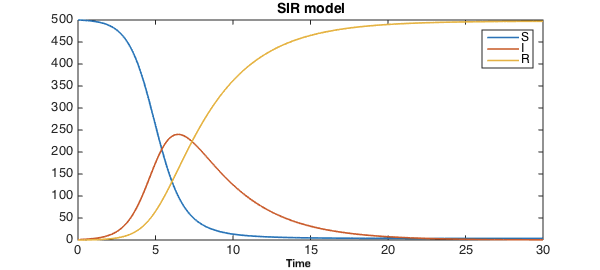}
	\smallskip
	
	\includegraphics[height=.25\textheight]{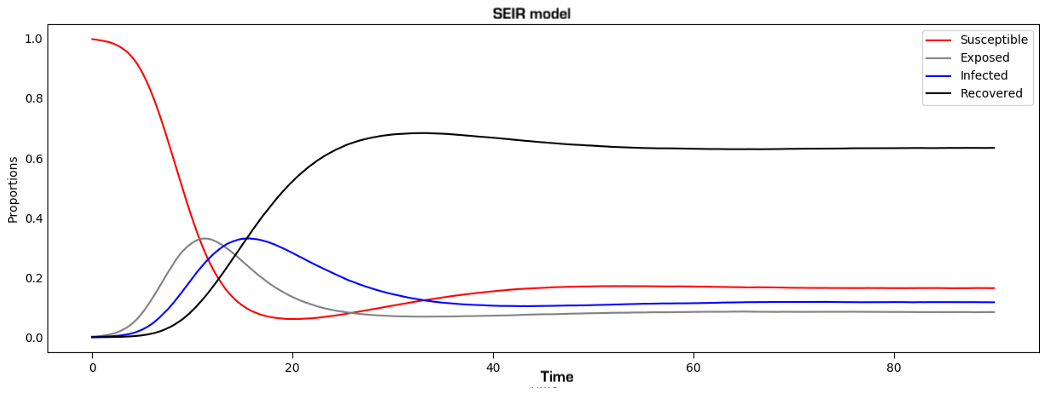}
	\caption{SIR (\cite{SIR-Model}) and SEIR (\cite{SEIR-Model}) models of Covid rates}
	\label{SIRandSEIR}
	\end{figure}

These models are given by differential equations and involve statistical functions that evolve through time. An important property for using particular differential equations to  study quantitative phenomena is if their solutions  preserve probabilities. Fortunately, many  physical phenomena  are modeled by collections of equations that have this property. Examples include  \emph{Hamiltonian Systems} of differential equations.

\noindent\textbf{The basic questions}
What is the predictive behavior of models derived from differential equations?  Can you tell whether two models predict the same or different behavior? Will they predict random behavior? Will the behavior be completely predictable from an observation at time zero?
\smallskip

\noindent\textbf{The historical projects} 
The projects proposed by von Neumann in the 1930's (\cite{vn}) and Smale in the 1960's (\cite{smale}) were to \emph{classify} the statistical and qualitative behavior of the respective systems.  The idea was to have a way of cataloguing or indexing all of the systems in a given class so that two systems that behave the same way were indexed identically.
\smallskip

\noindent\textbf{What does \emph{Same Behavior} mean?} Two transformations $S$ and $T$ have the same behavior if they are \emph{conjugate}.  The notion of conjugacy is a generalization of a change of coordinates whose precise meaning depends on whether the context is qualitative or quantitative.  The simplest example is a transformation stated in $\mathbb R^2$ in terms of $(x,y)$-coordinates and the coordinates get changed to polar $(r,\theta)$, and can be changed back again. \emph{Conjugacy} allows more general changes, but with a similar idea.

\emph{Classification} is understood to mean finding a way of discerning whether transformations given in different frameworks are conjugate. Informally two transformations are conjugate just in case they exhibit exactly the same behavior in the appropriate sense.

\emph{Invariants:} The usual approach to classification is to assign \emph{invariants} to transformations that are the same for all elements of a conjugacy equivalence class. A classical example is the Koopman operator. In the quantitative case, there is a Hilbert space $L^2(X)$ associated to the setting for the transformations and a transformation that preserves statistics is associated to a unitary operator $\mathbb K:L^2(X)\to L^2(X)$.  This operators has \emph{discrete spectrum} if there is a basis for $L^2(X)$ consisting of eigenvectors.
\smallskip

\noindent\textbf{The New Results}
The new results are that in many cases  the classification program is impossible in a rigorous sense in both types of 
systems. For the quantitative systems this is joint work with B. Weiss (\cite{anticlass}) and for the qualitative systems this 
is joint work with A. Gorodetski (\cite{FGarxiv}). In both cases it is shown that there is no classification possible using only  inherently countable information. 
The results are stated rigorously in section \ref{RM}.

\section{Dynamical Systems formalism}
 Systems whose states that evolve through time are called \emph{dynamical systems}. Let us state this 
 formally.  Let $X$ be a ``nice'' space.  In this context $X$ will be a manifold, usually compact, or $X$ will 
 be a probability space. 
 
 An $\mathbb R$-action or a \emph{flow} on $X$ is a function $T:\mathbb R\times X\to X$ such that for 
 all 
 $x\in X, s, t\in \mathbb R$:
 	\begin{equation}\label{action}
	T(s+t, x)=T(s, T(t,x)).
	\end{equation}
Conceptually: $t$ is the time variable. The equation is saying that the result of starting at $x$ and going 
to where $x$ is at time $s+t$ is that same as starting at $x$ and first going to where $x$ is at time $t$ to 
get a point $y=T(t,x)$ and then starting at $y$ and going to $T(s, y)$.   

\smallskip

\noindent\textbf{{$\mathbf{\mathbb Z}$} vs. $\mathbf{\mathbb R}$:} Modeling time continuously uses $\mathbb R$-flows. Pointing to work of Poincar\'{e} (\cite{Poinc}) and Birkhoff (\cite{19}) Smale  states  ``the same phenomena and problems of the qualitative theory of ordinary differential equations are present in their simplest form in the [\emph{single}] diffeomorphism problem.  Having first found theorems in the diffeomorphisms case, it is usually a secondary task to translate the results back to the" [\emph{continuous time framework ...}] (\cite{smale1967})
In other words one should study single transformations or equivalently $\poZ$-actions.   To connect with $\mathbb R$-flows, consider the single transformation $T$ with $T(x)=T(1,x)$. Then $T^n(x)=T(n,x)$. For the rest of the paper we discuss single transformations.

\smallskip

\noindent\textbf{What does \emph{ergodic} mean?}
For empirical systems an important property is that the average of repeated sample values of a function $f$ yields the average 
value  of $f$ over the whole domain under discussion:

	\begin{equation}\label{ergtheorem}
	\lim_{N\to \infty} \frac{1}{N}\sum_{n=0}^Nf(T^nx)= \int f(x)d\mu(x)
	\end{equation}
where $\mu$ is a probability measure that is invariant under $T$. 

Equation \ref{ergtheorem} is clearly implausible if a non-trivial part $I$ of the system evolves completely isolated from the rest of the system. Then sampling in $I$ would lead to global conclusions about the behavior of $X\setminus I$. In the context of Covid:  if country $A$ and country $B$ are completely isolated from each other with no possible sequence of transmission between $A$ and $B$, then sampling the rates of covid in $A$ would only give information about the rates of covid in $A$, not about the average rates of covid in both countries combined. The mathematical term for a system  $(T,\mu)$ not having
 non-trivial isolated invariant sets is that it is \emph{ergodic}.
For ergodic systems equation \ref{ergtheorem} is true.

\section{The quantitative study} The examples  of Lagrangians and SIRS models are well-understood compared to the general problem, which we now discuss in the quantitative case.  To frame the issue we consider two diametrically opposite special cases: the Bernoulli Shift and translations on compact groups.
\smallskip

\noindent\textbf{Completely Random Behavior} Let $H=\{1, 2, 3, 4, 5, 6\}$ and give each $n\in H$ value $1/6$ so $H$ can be viewed as a 6-sided die.  Let $X$ be the collection of functions $f:\poZ\to H$, and view $X$ as rolling the die infinitely far into the past and the future. $X$ can be made into a measure space by setting the probability of a given sequence of rolls $\la n_0, n_1, \dots n_k\ra$ to be $(1/6)^k$. The transformation defined on $X$ is the shift function:
	\[sh(f)(k)=f(k+1).\]
This models rolling the die once:  the current time-0 roll becomes the time $-1$-roll. When you roll again it becomes the time $-2$-roll. The system $(X, sh)$ is a model of the  completely random behavior of rolling dice. 

By generalizing this to an $M$ faced die,  having non-zero, possibly uneven, probabilities that sum to one, we 
get the notion of a Bernoulli shift. These systems are ergodic. They seem easy to understand, but 
Meshalkin (\cite{meshalkin}) proved that the system determined by a 5-headed die with probabilities 
$\{1/2, 1/8, 1/8, 1/8, 1/8\}$ is conjugate to the 4-headed die with probabilities $\{1/4, 1/4, 1/4, 1/4\}$.  
This makes the problem of determining conjugacy harder to understand.

Kolmogorov and Sinai adapted Shannon's information theoretic idea of \emph{entropy} to understanding signals. Viewing the rolls of the dice as a discrete signal one sees that the entropy of a $k$-headed Bernoulli shift with probabilities $\{p_1, p_2, \dots p_k\}$  (with $\sum_ip_i=1$) is
	\[h=-\sum_{i=1}^k log(p_i).\]
\smallskip

The problem of classifying  Bernoulli shifts was solved by Ornstein \cite{orn}:
\medskip

\noindent\textbf{Theorem}
\emph{Two Bernoulli shifts are conjugate if and only if they have the same entropy.}
\smallskip

\noindent In the language of this paper, the entropy function computes \emph{complete numerical invariants}.

\medskip

\noindent\textbf{Completely Predictable Behavior} Let $G$ be a compact topological group, such as the unit circle or the torus. Then $G$ carries a translation invariant probability measure, called \emph{Haar measure}. Suppose that $a$ is an element of $G$ whose powers are dense in $G$.  Then $G$ must be abelian and the map sending $g$ to $g+a$ is ergodic.

 Halmos and von-Neumann (\cite{HvN}) showed:
 \medskip
 
 	\noindent\textbf{Theorem}
	\emph{Let $X$ be a measure space and $T:X\to X$. Then the following are equivalent:
		\begin{enumerate}
		\item $T$ is conjugate to an ergodic shift on a compact group $G$,
		\item The Koopman operator associated with $T$ has discrete spectrum.
		\end{enumerate}
	If either of these conditions hold then  both $G$ and  the transformation $T$ are determined by the countable collection of eigenvalues of the Koopman operator} 
\medskip

\noindent\textbf{How do these examples differ?} 
\begin{itemize}
	\item In the first example the entropy function $h$ is a lower semi-continuous function that 
	assigns a real number to each Bernoulli shift. Two shifts are conjugate just in case they are 
	assigned the same number.
	
	\item In the example of translations on a compact groups, there is a function $E$ that 
	computes the countable group of eigenvalues $\la\lambda_n:n\in\nn\ra$ of the unitary operator
	 and enumerates them in a specific order. Given conjugate transformations $S$, $T$, 
	 $E(S)$ and $E(T)$ enumerate the same collection of  eigenvalues  but possible in a 
	 different order. 
\end{itemize}
Thus both examples have flaws: in the first example, the given complete numerical invariant 
is not recursively computable.  In the second example the invariant is a countable \emph{set}, but the set doesn't have a standard way of enumerating it.  So the invariant determining the conjugacy class is not even a number! It turns out that the situation for general measure preserving transformations, even if they are very nice diffeomorphisms, is much worse.

The main question  is:
	\begin{quotation}
	\noindent How effectively can we tell whether two measure preserving diffeomorphisms are conjugate?
	\end{quotation}
If we can't even tell diffeomorphisms apart in general, we certainly can't hope to isolate the behavior of a particular diffeomorphism except in special cases.  The examples of the 2-body problem and the SIRS models given above are those special cases.  But even for general Hamiltonians, or the $n$-body problem it is not known.

\section{Is the problem even possible to solve?}  Is it \emph{feasible} to tell two transformations apart? There is a general theory of feasibility, that we now outline.
\medskip

\noindent Measures of feasibility in decreasing strength include the properties:
 	\begin{enumerate}
	\item having an algebraic solution in closed form,
	\item having an algorithm solving a problem that runs in a reasonable time,
	\item having \emph{any} algorithm that solves a problem or allows an effective classification,
	\item having a solution that uses only inherently countable information,
	\item solutions that use uncountable Axiom of Choice
	\end{enumerate} 
	
	The \textbf{3-body problem} is an  example of an elementary problems that fails the first property. However in 
	1907 Sundman  did show it satisfies property 4 (\cite{sundman}).  Other examples failing the first item include
	 computing $\int e^{-x^2}dx$ and from elementary Galois Theory.
\smallskip  

Many problems fail the second property.  For example, one important idea behind \textbf{Cryptocurrency} is that even though there is an  algorithm for finding the key, that the algorithm can't run in reasonable time. 
\smallskip  

The classical \textbf{word problem} in group theory is an example of a problem that has no algorithm to solve, hence fails property 3. An appeal to Church's Thesis allows one to restate the third property as saying that a problem is 
solvable using \emph{inherently finite} information. 
\smallskip 

Property 4 is a very significant weakening of property 3, since it allows infinite information, as long as it is inherently countable. Thus saying that there is no solution that uses inherently countable information is a very strong negative statement.  It is property 4 that we focus on.
\smallskip

The uncountable Axiom of Choice is viewed as controversial by some people, precisely because 
it has \emph{no} element of feasibility. Without the uncountable Axiom of Choice \textbf{it is impossible 
to construct a non-measurable set} (\cite{solovay}).
\smallskip 

While not really a measure of effectiveness, there is an intermediate notion between property 4 and the uncountable axiom of choice: the notion of being analytic. A set is \emph{analytic} if it is the continuous image of a closed subset of a Polish space. Suslin (\cite{Sus}), correcting a famous mistake of Lebesgue (\cite{Leb}), showed that this is more general than property 4.

There are also degrees of ineffectiveness: there is a partial ordering $\le_\mcb$ such that if $A\le_\mcb B$ then $B$ is at least as ineffective as $A$. Among the analytic sets there  are 
$\le_\mcb$-maximal sets. These  can be viewed as maximally complex analytic sets. There are not computable using inherently countable information. 

\section{How feasible is classifying the quantitative and qualitative behavior?}\label{RM}
Both the collection of infinitely differentiable diffeomorphisms of a fixed manifold and the collection of infinitely differentiable measure preserving diffeomorphisms of a manifold carrying a standard measure nice topologies.  So the questions of effectiveness outlined in the previous section can be asked in these contexts.
We address two questions.  Let $X$ be either of these two spaces and let $E$ be either the equivalence relation of conjugate-by-a-homeomorphism (in the qualitative case) 
or conjugate-by-a-measure preserving transformation (in the quantitative case). 
\medskip

\noindent{\emph{Question 1}} Is there a function $f:X\to \mathbb R$ that uses only countable information
	such that $x_1Ex_2$ if and only if $f(x_1)=f(x_2)$?  This can be re-stated as asking whether it is possible to have complete numerical invariants for the equivalence relation $E$ that can be computed using countable resources

	\medskip

\noindent{\emph{Question 2}} Can the question of whether $x_1 E x_2$ be settled at all with countable information?

\medskip

\begin{center}
\textbf{The answers to both questions are negative in both the quantitative and the qualitative setting.} 
\end{center}

\noindent\textbf{The quantitative setting.}
 The first result along these lines in the context of dynamical systems is due to the author and A. Louveau (\cite{descview}). 
 \smallskip
 
\noindent{\textbf{Theorem}}\emph{(Foreman, Louveau) Let $\mathcal G$ be the collection of ergodic translations of 
compact groups.  Then there is no method of computing complete numerical invariants using inherently countable information.}
\smallskip

It follows that even in the context of the Halmos-von Neumann theorem, one can't measure complexity with numbers. This contrasts with Ornstein's theorem where entropy is a complete numerical invariant.

But what about the general situation? Can you distinguish between ergodic diffeomorphisms using inherently countable methods? The answer is no:
\smallskip

\noindent{\textbf{Theorem}}\emph{(Foreman, Weiss)(\cite{anticlass}) Let $E$ be the equivalence relation of conjugacy-by-a-measure-preserving-transformation on diffeomorphisms of the 2-torus. Then there is no inherently countable way of determining whether two diffeomorphisms are conjugate. In fact $E$ is maximally complex.}

\begin{figure}[!h]
	\centering
	\includegraphics[height=.25\textheight]{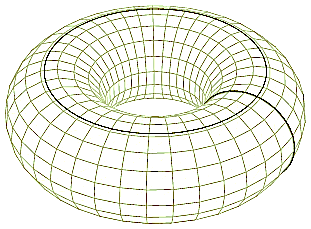}
	\caption{Image of the 2-torus (\cite{2torus})} 
	\label{2torus}
	\end{figure}

\noindent\textbf{The qualitative setting} The situation for diffeomorphisms with the equivalence relation of conjugacy-by-homeomorphisms is more subtle. Large classes have been classified, especially in lower dimensions. Work of Poincar\'{e} \cite{Poinc} and later Denjoy \cite{denjoy} effectively classified aperiodic diffeomorphisms on one dimensional compact manifolds (the circle). In 2 dimensions there are very large classes--for example the Morse-Smale diffeomorphisms are examples.

However, the general situation is quite similar to the quantitative study.
\smallskip

\noindent{\textbf{Theorem}}\emph{(Foreman, Gorodetsky)(\cite{FGarxiv}) Let $M$ be a manifold of dimension at least 2.  Then there is no method of computing complete numerical invariants for the qualitative behavior of diffeomorphisms using inherently countable information.}
\smallskip

The question about the equivalence relation has a similar answer.

\noindent{\textbf{Theorem}}\emph{(Foreman, Gorodetsky)(\cite{FGarxiv}) Let $M$ be a manifold of dimension at least 5.  Then the equivalence relation of topological conjugacy is not Borel.  In fact $E$ is maximally complex.}

\noindent Quite recent work has led to results by the same two authors that seem to show the analogous results in all dimensions.

\noindent\textbf{But the objects must be very strange!} One might think that an impossible--to--classify collection of objects 
must be on exotic spaces and consist of extreme points.  However this is not true. In the quantitative case, they are 
completely smooth ($C^\infty$) diffeomorphisms of the 2-torus. (See figure \ref{2torus}.)

Moreover it is impossible to classify diffeomorphisms on  \emph{any} 2-dimensional manifold with a free action of the circle group.  
The behavior described in this section is typical, rather than exotic.

\medskip
In the qualitative case, inside the space of diffeomorphisms of 2-manifolds, there are open sets that do have good classifications, such as the Morse-Smale diffeomorphisms (\cite{Peix}).  But these seem to be islands of order in a sea of chaos. 

\medskip
\noindent\textbf{Has all the work been done?} By no means.  Large classes of measure preserving diffeomorphisms, including the weakly mixing diffeomorphisms, are thought to be unclassifiable. While there are tentative results in this direction, they have not been verified.

%
%
%

\bibliography{arch.bib}
 
 \bibliographystyle{plain}

 \end{document}